\documentstyle[12pt,twoside]{amsart}

\title{The global indices of log Calabi-Yau varieties\\
--A supplement to Fujino's paper: The indices of log canonical 
singularities--
}
\author{Shihoko ISHII}

\address{Shihoko Ishii: Department of Mathematics,
         Tokyo Institute of 
         Technology,
         Oh-Okayama, Meguro, Tokyo, Japan}
\newcommand{\bN}{{\Bbb N}}
\newcommand{\bC}{{\Bbb C}}
\newcommand{\bQ}{{\Bbb Q}}
\newcommand{\bP}{{\Bbb P}}
\newcommand{\Ind}{{\operatorname{Ind}}}
\newcommand{\im}{{\operatorname{Im}}}

\newtheorem{thm}{Theorem}[section]

\newtheorem{lem}[thm]{Lemma}

\newtheorem{prop}[thm]{Proposition}

\theoremstyle{definition}
\newtheorem{defn}[thm]{Definition}

\newtheorem{say}[thm]{}
\newtheorem{exmp}[thm]{Example}

\newtheorem{rem}[thm]{Remark}

\theoremstyle{remark}

\begin{document}
\maketitle
{\footnote[1]{Partially supported by the Grant-in-Aid for Scientific
Research, the Ministry of Education, Japan.}}

\pagestyle{myheadings}
\markboth{\hfill SHIHOKO ISHII\hfill}{\hfill GLOBAL INDICES OF 
LOG CALABI-YAU \hfill}
\begin{abstract}
  This paper gives the all possible global indices of
  log Calabi-Yau 3-folds with standard coefficients on the boundaries
  and having lc, non-klt singularities.
  This follows easily from the discussion in the paper:
  The indices of log canonical singularities by Fujino.
\end{abstract}

\section{Introduction}
  In this paper,  we study a log pair \( (X, B_{X}) \)
  with a normal projective variety \( X \) defined over \( \bC \) and a 
  boundary \( B_{X} \)
  of standard coefficients (i.e., \( B_{X}=\sum b_{i}B_{i} \), where
  \( b_{i}=1 \) or \( 1-1/m \) for \( m\in \bN \)).
  A pair \( (X,B_{X}) \) is called a log Calabi-Yau variety if 
  it has  lc singularities and
  \( K_{X}+B_{X}\equiv 0 \). 
  For a log Calabi-Yau variety \( (X,B_{X}) \) assume that there exists \( 
  r\in \bN \) such that \( r(K_{X}+B_{X})\sim 0 \).
  (For \( \dim X\leq 3 \) this holds true for every log Calabi-Yau 
  variety,
  by the abundance theorem ({\cite[11.1.3]{utah}}    and \cite{k-m-m}).
  We define the global index \( \Ind(X,B_{X}) \) by the minimum of such \( r \).

  It is well known that a non-singular surface \( X \) with \( 
  K_{X}\equiv 0 \)   has \( \Ind(X,0)=1,2,3,4,6 \).
  Blache \cite{blache} proved 
  that a normal surface  \( X \) with \( 
  K_{X}\equiv 0 \) and having lc non-klt singularity
  has also \( \Ind(X,0)=1,2,3,4,6 \). 
  This is generalized into the case that  log Calabi-Yau surface \( (X,B_{X}) \)
  has lc and non-klt singularities in {\cite[2.3]{shokurov}}.

  In this paper we prove the following:
  
\begin{thm}
\label{3-fold}
  Let \( (X, B_{X}) \) be a  log Calabi-Yau 3-fold with lc
  non-klt singularities. 
  Then \( r\in \bN \) can be the global index \( \Ind(X, B_{X}) \),
  if and only if \( \varphi(r)\leq 20 \) and \( r\neq 60 \),
  where \( \varphi \) is the Euler function.
  In particular the global index is bounded.
\end{thm} 

This theorem is a corollary of the following:

\begin{thm}
\label{main}
  Assume the Abundance Theorem and \( G \)-equivariant log Minimal 
  Model Program for dimension \( \leq n \),
  where \( G \) is a finite group.
  Let \( (X, B_{X}) \) be an \( n \)-dimensional log Calabi-Yau 
  variety with non-klt singularities.
  If the conjectures \( (F'_{j}) \) and \( (F_{l}) \) in \cite{fujino}
   hold true for \( j= n-1 \), \( l\leq n-2 \),
   then the global index \( \Ind(X,B_{X}) \) is bounded.
\end{thm} 

  The author would like to express her gratitude to Prof. Yuri Prokhorov
  for calling her  attention to this problem and offering useful comments
  and stimulating discussions.
  She also express her gratitude to Dr. Osamu Fujino and Prof. Pierre Milman
   for giving useful 
  information about their papers.
  Dr. Osamu Fujino also pointed out a mistake of the preliminary 
  version of this paper, for which she is grateful to him.  

\section{The global indices}
\begin{say}
\label{assumption}
  Throughout this paper,
  we use the notation and the terminologies 
  in \cite{fujino}.
  We assume the Abundance Theorem and the \( G \)-equivariant log 
  Minimal Model Program (as is well known, these hold for dimension \( 
  \leq 3 \) by {\cite[11.1.3]{utah}, \cite{k-m-m} and 
  \cite[2.21]{k-m}}).
\end{say}  
\begin{say}
  Let \( (X, B_{X}) \) be an \( n \)-dimensional log Calabi-Yau 
  variety. 
  Since we assume the Abundance Theorem,
  there exists \( r\in \bN \) such that \( r(K_{X}+B_{X})\sim 0 \).
  Let \( \pi:(Y,B)\to (X,B_{X}) \) be the  index 1 cover with 
  \[  K_{Y}+B=\pi^*(K_{X}+B_{X}) .    \] 
  Here the  index 1 cover is constructed as follows:
  let  \( r=\Ind(X, B_{X}) \),
  then there exists a rational function \( \varphi \) on \( X \) such 
  that \( r(K_{X}+B_{X})=\operatorname{div}(\varphi) \);
  take the integral closure \( Y \) in \( 
  K(X)(\sqrt[r]{\varphi}) \). 
  Note that \( K_{Y}+B\sim 0 \), that  \( B=\pi^*(\lfloor B_{X}\rfloor) \)
  is a reduced divisor and that \( \pi \) ramifies only over the components of \( 
  B_{X} \) whose coefficients are \( <1 \), 
  as the coefficients of \( B_{X} \) are standard.
  Since \( K_{X}+B_{X} \) is lc (resp. klt) if and only if \( K_{Y}+B \)
  is lc (resp. klt),
  \( (Y, B) \) is  log Calabi-Yau of global index 1. 
  Therefore we obtain that every log Calabi-Yau variety
  \( (X,B_{X}) \) is the quotient of a log Calabi-Yau variety of 
  global 
  index 1 by the action of a finite cyclic group.
\end{say}

\begin{say}  
  Let \( G \) be the cyclic group acting on a log Calabi-Yau variety 
  \( (Y,B) \) of global index 1.
  Since \( G \) acts on \( \Gamma(Y, K_{Y}+B)=\bC \),
  there is a corresponding representation \( \rho : G\to GL(\Gamma(Y, 
  K_{Y}+B))=\bC^* \).  
\end{say}

\begin{lem}
\label{key}
  Under the notation above,
  Let \( (X, B_{X}) \) be the quotient \( (Y, B)/G \) by \( G \).
  Then 
  \[ \Ind(X, B_{X})=| \im \rho | . \]
\end{lem}  
      
\begin{pf} 
  For a generator \( \theta\in \Gamma(Y, K_{Y}+B) \),
  \( \theta^{|\im \rho|} \) is \( G \)-invariant,
  therefore \( \Gamma(X, | \im\rho | (K_{X}+B_{X}))\neq 0 \),
  which yields \( \Ind(X, B_{X})\leq | \im \rho |   \).
  Conversely,
  for a generator \( \eta \in \Gamma(X, \Ind(X,B_{X})(K_{X}+B_{X})) 
  \),
  \( \pi^*\eta\in \Gamma(Y, \Ind(X,B_{X})(K_{Y}+B))  \) is \( G 
  \)-invariant.
  If we write \( \pi^*\eta=a\theta^{\Ind(X,B_{X})} \) \( (a\in \bC) \),
  for a generator \( g\in G \), \( (a\theta^{\Ind(X,B_{X})})^g=
  a\epsilon ^{\Ind(X,B_{X})}\theta^{\Ind(X,B_{X})}=
  a\theta^{\Ind(X,B_{X})} \), where \( \epsilon \) is a primitive 
  \( |\im \rho |  \)-th root of unity.
  Hence,   \( \Ind(X, B_{X})  \geq |\im \rho|   \). 
\end{pf}
    
\begin{say}
\label{nonklt}
  Now we are going to study lc and non-klt log Calabi-Yau 
  varieties.
  Let \( (Y, B) \) be an \( n \)-dimensional log Calabi-Yau variety of 
  global index 1  with lc and non-klt singularities.
  Assume that a cyclic group \( G \) acts on \( (Y, B) \).
  Then we have a  projective \( G \)-equivariant log resolution 
  \( \varphi: {\tilde Y} \to Y \) of \( (Y,B) \).
  Indeed,  
  let \( \varphi': {\tilde Y}' \to Y \) be the canonical resolution of \( 
  (Y,B) \) constructed in \cite{b-m},
  then \( \varphi' \) is projective and \( {\varphi'}^{-1}(B)\cup \)(the 
  exceptional set)  is  normal crossing divisor.
  By the blow up at a suitable \( G \)-invariant center,
  we obtain the divisor with simple normal crossings.
  Define the subboundary \( F \) on \( {\tilde Y} \) by 
  \( K_{\tilde Y}+F=\varphi^*(K_{Y}+B) \).
  Run \( G \)-equivariant log MMP for \( K_{\tilde Y}+F^B \) over \( Y \)
  (The notation \( F^B \) is in \cite[1.5]{fujino} and \( F^B=F^c \) 
  in our case).
  Then we obtain \( G\bQ \)-factorial dlt pair \( f: (Y',B') \to 
  (Y,B) \) over \( (Y,B) \).
  Since \( K_{Y'}+B' \) is \( f \)-nef and  \( (Y,B) \) is lc,
  we obtain that \( K_{Y'}+B'=f^*(K_{Y}+B)\sim 0 \).
  By \cite[2.4]{fujino}, \( B' \) has at most two connected 
  components.          
\end{say}
 
\begin{defn}[for the local version, see {\cite[4.12]{fujino}}]
  Let \( (Y, B) \) and \( ({\tilde Y},F) \) be as in \ref{nonklt}.
  We define
  \[ \mu=\mu(Y,B):=\min\{\dim W \mid W\in CLC({\tilde Y}, F)\}. \]
\end{defn}

Note that in case \( B' \) is connected, then \( 0\leq \mu\ \leq n-1 \)
and in case \( B' \) has two connected components, then \( \mu = n-1 \).
  
{\bf Case 1} (\( B' \) is connected)

  There exist a \( G \)-isomorphism \( \Gamma(Y,K_{Y}+B)\simeq 
  \Gamma(Y',K_{Y'}+B') \)  and an exact sequence:
  \[ 0=\Gamma(Y', K_{Y'})\to \Gamma(Y',K_{Y'}+B')
  \to \Gamma(B', (K_{Y'}+B')|_{B'})=\bC, \]
  where the last term is isomorphic to \( \Gamma(B', K_{B'}) \),
  as \( K_{Y'}+B' \) is a Cartier divisor.
  Therefore, we have only to check the action of \( G \) on 
  \( \Gamma(B', K_{B'}) \).

\begin{prop}[for the local case, see {\cite[4.11]{fujino}}]
\label{connected0}
  If there exists a non-zero admissible section in 
  \( \Gamma(B', m_{0}K_{B'}) \), then \( G \) acts on 
  \( \Gamma(B', m_{0}K_{B'}) \) 
  trivially.
\end{prop}    

\begin{pf}
  The proof is the same as that of {\cite[4.11]{fujino}}.
  We have only to note that \( B'=E=E^c \) in our case.
\end{pf}

\begin{prop}[for the local case, see {\cite[4.14]{fujino}}]
\label{connected1}
  Assume that \( \mu(Y,B)\leq n-2 \). 
  Then there exists a non-zero admissible section 
  \( s\in \Gamma(B', m_{0}K_{B'}) \) with \( m_{0}\in D_{\mu} \).
  In particular, \( s \) is \( G \)-invariant.
  Thus, \( \Ind ((Y,B)/G)\in I_{\mu} \).
\end{prop}
  
\begin{pf}
  The proof is the same as that of  {\cite[4.14]{fujino}}. 
  Again \( B'=E=E^c \).
\end{pf}

\begin{prop}
\label{connected2}
  Assume that \( B' \) is connected and \( \mu(Y,B)=n-1 \).
  Then \( \Ind(Y,B)/G\in I'_{n-1} \).
\end{prop}

\begin{pf}
  In this case,
  \( B' \) is irreducible, therefore \( (Y', B') \) is plt.
  Then, by Adjunction  {\cite[17.6]{utah}}, \( B' \) is klt and \( K_{B'}\sim 0 \).
  Now apply \ref{key}.
\end{pf}

{\bf Case 2} (\( B' \) has two connected components).

  Note that \( B' \) is the disjoint union of two irredicible
  components, therefore \( (Y', B') \) is plt (see 
  \cite[2.4]{fujino}).
  Run  \( G \)-equivariant log MMP for \( K+B'-\epsilon B' \),
  then we obtain a  \( G \)-equivariant contraction 
  \( p:Y''\to Z \)  of an extremal face for \( K+B''-\epsilon B'' \)
  to a lower dimensional variety \( Z \),
  where \( B''= B''_{1}\amalg B''_{2} \) is the  divisor on \( Y'' \)
  corresponding
  to \( B' \).
  Here \( \dim Z=n-1 \),
  because \( h^{n-1}(Y', {\cal O}_{Y'})=h^1(Y', K_{Y'})\neq 0 \).    
  We also obtain that \( B''_{i} \)'s are generic sections of \( p \).
  Since \( (Y'',B'') \) is plt and \( K_{Y''}+B''\sim 0 \),
  each \( B''_{i} \) has canonical singularities and \( K_{B''_{i}}\sim 
  0 \) by {\cite[17.6]{utah}}.
  Then the birational image \( Z \) has \( K_{Z}\sim 0 \),
  and therefore it has canonical singularities.
   Since the group \( G=\langle g \rangle \) acts on \( B'' \),
  the subgroup \( H:=\langle g^2 \rangle \)   acts on each \( B''_{i} \) \( (i=1,2) \).
     Consider the exact sequence:
   \[ 0=\Gamma(Y'',K_{Y''}+B''_{2})\to \Gamma(Y'', 
   K_{Y''}+B'')\stackrel{\alpha}\longrightarrow 
   \Gamma(B''_{1},K_{B''_{1}}), \]
   where \( \alpha \) is an \( H \)-equivariant isomorphism.
  On the other hand, 
  the homomorphism 
  \( \Gamma(B''_{1}, K_{B''_{1}}) \to \Gamma(Z, K_{Z}) \)
  induced from \( p|_{B''_{1}} \)  
  is also an \( H \)-equivariant isomorphism.
  Hence, for two representations \( \rho:G \to GL(\Gamma(Z, K_{Z})) \)
  and \( \rho':G \to GL(\Gamma(Y'', K_{Y''}+B'')) \),
  we obtain the equality \( |\rho(H)|=|\rho'(H)| \).
  Note that, for any representation \( \lambda:G \to {\Bbb C}^* \),  
  \( \lambda(H)=\lambda(G) \) if and only if \( |\lambda(G)| \) is an 
  odd number.
  If we denote \( |\rho(G)| \) by \( r \),
  then \( r\in I'_{n-1} \), and either:
  (1) \( |\rho'(G)|=r \) or (2) \( |\rho'(G)|=2r \) and \( r \) is odd 
  or (3) \( |\rho'(G)|=r/2 \) and \( r/2 \) is odd.
  By defining \( I''_{k}:=I'_{k}\cup \{2r \mid r\in I'_{k} \ {\rm 
  odd}\} \cup\{r/2 \mid r\in I'_{k}, r/2 \ {\rm odd}\} \),
  we obtain:
  
  \begin{prop}
\label{non-connected}
  Assume \( B' \) has two connected components.
  Then \( \Ind((Y,B)/G)\in I''_{n-1} \).
\end{prop}

  By  \ref{connected1}, \ref{connected2} and 
  \ref{non-connected},
  we obtain Theorem \ref{main}.
  In particular, for the 3-dimensional case \( G \)-equivariant log MMP,
  the Abundance Theorem and  \( (F'_{j}) \), \( (F_{l}) \) \( (j=2, 
  l\leq 1) \) hold.
  Here note that \( I_{0}=\{1,2\} \), \( I_{1}=\{1,2,3,4,6\} \) and
  \( I'_{2}=\{r\in {\Bbb N} \mid \varphi (r)\leq 20, r\neq 60\} \) by 
  \cite{nikulin} and \cite{m-o}.
  By the list of the values of \( I'_{2} \) in {\cite[Table 1]{m-o}},
  we can check that \( I''_{2}=I'_{2} \).
  Therefore  we obtain the necessary condition of 
  the global index
  \( \Ind(X, B_{X}) \)  in Theorem \ref{3-fold}.  

  The following shows that it is the sufficient condition of the 
  global index:
  
\begin{exmp}
  Let \( r \) be a positive integer that satisfies \( \varphi(r) \leq 20 \)
  and \( r \neq 60 \).
  Then by \cite{kondo} and \cite{m-o},
  there exists a \( K3 \)-surface \( S \) with an action \( G \) of order \( r \)
  and \( r = |\im \rho| \).
  Let \( Y = S\times \bP^1\) and \( B = S\times \{0\}+S \times 
  \{\infty\} \).  
  Let \( G \) act on \( Y \) by trivial action on \( \bP^1 \) and
  the action above on \( S \).
  Let \( (X, B_{X}) \) be the quotient of \( (Y, B) \) by \( G \)
  with \( K_{Y}+B=\pi^*(K_{X}+B_{X}) \).
  Then  \( (X, B_{X}) \) is a log Calabi-Yau 3-fold with global index \( 
  r \).     
\end{exmp}

\begin{rem}
  We can also prove Theorem \ref{3-fold} by using  \cite{ishii} 
  instead of \cite{fujino}.
  Indeed, we used \cite{fujino} only for propositions \ref{connected0}
  and \ref{connected1}.
  For the 3-dimensional case, these propositions can be replaced by
  the discussion on the order of the action of \( G \) on 
  \( H^2(F^B, {\cal O}_{F^B}) \) for type \( (0, 0) \) and \( (0, 1) \).
  Theorems {\cite[4.5]{ishii}}  and {\cite[4.12]{ishii}}  give the same results as in 
  \ref{connected1}.      
\end{rem}

\begin{rem}
  Osamu Fujino informed the arthor that the boundedness of the indices 
  of log Calabi-Yau 3-folds also follows from \cite[4.17]{fujino} and 
  the proof of  \cite[4.14]{fujino}.
  By this proof we obtain the index in \( I_{2} \) instead of \( I'_{2} \).   
\end{rem} 

\begin{rem}
  If we assume \( (F'_{n}) \),
  then it is clear that \( n \)-dimensional klt log Calabi-Yau variety
  has the global index \( r \in I'_{n} \) by Lemma \ref{key}.
  Therefore klt log Calabi-Yau surface has the global index \( r \)
  such that \( \varphi(r)\leq 20 \) and \( r \neq 60 \).
  
   For a klt log Calabi-Yau 3-fold with \( B_{X}=0 \),
    the global index satisfies the same condition as 
  above {\cite[Corollary 5]{m-o}}.
  
\end{rem}

\makeatletter \renewcommand{\@biblabel}[1]{\hfill#1.}\makeatother

\end{document}